\documentclass{aptpub}
\usepackage{graphicx}
\usepackage{amsmath}
\usepackage{amsfonts}
\usepackage{amssymb}
\usepackage{setspace}
\graphicspath{{./}{./figures/}{../subpoisson/figures/}{../R/}}
\usepackage[pdftex,colorlinks,%
           bookmarks=true,%
           pdftitle=On\ comparison\ of\ clustering\ properties\ of\ %
           point\ processes,%
           pdfauthor=B.\ Blaszczyszyn%
          \ -\  D.\ Yogeshwaran]{hyperref}
\hypersetup{
   bookmarksnumbered,
   pdfstartview={FitH},
   citecolor={blue},
   linkcolor={red},
   urlcolor=[rgb]{0,0.55,0},
   pdfpagemode={UseOutlines}
}

\newcommand{\be}{\begin{equation}}
\newcommand{\ee}{\end{equation}}
\newcommand{\bea}{\begin{eqnarray}}
\newcommand{\beas}{\begin{eqnarray*}}
\newcommand{\no}{\nonumber}
\newcommand{\eea}{\end{eqnarray}}
\newcommand{\eeas}{\end{eqnarray*}}

\newcommand{\qed}{\hfill $\Box$}
\newcommand{\EXP}[1]{\mathsf{E}\!\left(#1\right) }
\newcommand{\COV}[1]{\mathsf{Cov}\left( #1 \right)}

\newcommand{\remove}[1]{}

\newcommand{\pr}[1]{\mathsf{P}\left( #1 \right)}
\newcommand{\md}{\text{d}}

\newcounter{cnt1}

\newcounter{cnt3}
\newcommand{\blr}{\begin{list}{$($\roman{cnt1}$)$}
 {\usecounter{cnt1} \setlength{\topsep}{0pt}
 \setlength{\itemsep}{0pt}}}
\newcommand{\bla}{\begin{list}{$($\betaph{cnt2}$)$}
 {\usecounter{cnt2} \setlength{\topsep}{0pt}
 \setlength{\itemsep}{0pt}}}
\newcommand{\bln}{\begin{list}{$($\arabic{cnt3}$)$}
 {\usecounter{cnt3} \setlength{\topsep}{0pt}
 \setlength{\itemsep}{0pt}}}
\newcommand{\el}{\end{list}}

\def\rar{\rightarrow}

\newcommand{\al}{\alpha}

\newcommand{\lam}{\lambda}

\newcommand{\1}{{\bf 1}}

\def\mR{\mathbb{R}}
\def\mZ{\mathbb{Z}}
\def\mN{\mathbb{N}}

\def\mE{\mathbb{E}}
\def\mM{\mathbb{M}}

\newcommand{\cN}{{\mathcal N}}
\newcommand{\cX}{{\mathcal X}}

\newcommand{\ur}{{\underline r}}
\newcommand{\ovr}{{\overline r}}
\def\mL{\mathbb{L}}

\authornames{B. Blaszczyszyn , D. Yogeshwaran}
\shorttitle{On comparison of clustering properties
of  point processes}

\begin{document}


\title{On comparison of clustering properties
of  point processes}\\[2ex]
\authorone[INRIA/ENS]{Bart{\L}omiej  B{\L}aszczyszyn\vspace{-5ex}}
\addressone{23 av. d'Italie, 75214 Paris Cedex 13, FRANCE; email: Bartek.Blaszczyszyn@ens.fr}
\noindent\authortwo[Technion --- Israel Institute of Technology]{D. Yogeshwaran}
\addresstwo{Dept. of Electrical Engineering,
Technion --- Israel Institute of Technology,
Haifa 32000, ISRAEL; e-mail: yogesh@ee.technion.ac.il}
\vspace{-4ex}

\pdfbookmark[1]{Titlepage}{title} 
%


\begin{abstract}
In this paper, we propose a new comparison tool for spatial
homogeneity of point processes, based  on the joint examination of
void probabilities and factorial moment 
measures. We prove that  determinantal and permanental processes,
as well as, more generally, negatively and positively 
associated point processes are comparable in this sense to the
Poisson point process of the same mean measure.
We provide some motivating results on percolation and coverage processes and preview further ones on other stochastic geometric models
such as minimal spanning forests, Lilypond growth models, random simplicial complexes showing that the new tool is relevant for a systemic approach to the study of 
macroscopic properties of non-Poisson point processes.
This new comparison is also implied by the  directionally convex ($dcx$) 
ordering of point processes, which has already been shown to be relevant 
to comparison of spatial homogeneity of point processes. For this
latter ordering,  
using a notion of lattice perturbation, we provide
a large monotone spectrum of comparable point processes,
ranging from periodic grids to Cox processes, and encompassing Poisson
point process as well. They are intended to serve 
as a platform  for further theoretical and numerical 
studies of  clustering, 
as well as simple models  of random point patterns 
to be used in applications where neither complete regularity nor
the total independence property are realistic assumptions.
\end{abstract}

\setcounter{footnote}{0}
\renewcommand\thefootnote{\arabic{footnote}}

\section{Introduction}
\label{sec:intro}

Usual statistical approach to the study of  
clustering  in point processes (pp) 
consists in  the evaluation of  
Ripley's~$K$ function, pair-correlation function, or contact
distribution function (also called the empty space function).
However, such a  comparison  of local characteristics
seems a weak tool for the study of the impact of clustering on some macroscopic
properties of pp such as those required in
continuum percolation models. 
We are particularly motivated by heuristics indicating 
that pp exhibiting more  clustering should have 
larger critical radius for the percolation 
of its spherical-grain Boolean model 
than a spatially homogeneous pp.

It was observed in~\cite{snorder}, that the directionally convex ($dcx$)
order on pp implies the ordering of $K$ functions as well 
as  pair-correlation functions, in the sense that 
pp larger in the $dcx$ order have larger K functions and pair correlation functions,
while having the same mean number of points in any given set. 
Unfortunately, the examples from~\cite{snorder} 
are mostly only  some doubly-stochastic Poisson pp, which
are $dcx$ larger than Poisson pp (we call them
{\em super-Poisson} in this article).
In order to provide more examples of $dcx$ ordered pp, in particular
smaller than Poisson (we call them sub-Poisson), 
we study in this paper a notion of {\em perturbation of a pp}
consisting of  independent replication and translation of
points from some given, original pp.
A key observation is that such a perturbation  
is $dcx$  monotone with respect to the convex order
on the number of point replications.
In particular, perturbing a deterministic lattice in the above sense,  
one can obtain examples of  both sub- and super-Poisson pp,
with the Poisson pp itself obtained when 
the number of point replications has a Poisson distribution.
We believe these examples  can be useful
for modeling  of real phenomena  
for which neither lattice nor Poisson assumptions
can be justified. In this paper, we will also use them
to illustrate the aforementioned heuristic on 
the impact of clustering on the percolation of Boolean models.  

However, many examples of pp considered as clustering less or more
than the Poisson pp  of the same intensity
escape from the $dcx$ comparison; For example, 
{\em determinantal} and {\em permanental}
pp (cf.~\cite{Ben09}). In fact, despite  some structural similarities of
these pp to the perturbed lattices, we are able to
show for them $dcx$ order only on mutually disjoint  simultaneously 
observable sets, and not on all bounded Borel sets, required for the full $dcx$ order. 

The properties of {\em positive} and {\em negative
association} (cf~\cite{BurtonWaymire1985,pemantle00})
are also used to define classes of pp that, respectively,
cluster more or less than the completely independent (i.e., Poisson) pp.  
But it is not known if these properties imply or are implied by the $dcx$
ordering with respect to Poisson pp. Though one suspects many pp such as determinantal or
hard-core pp should be negatively associated, it is not known if they actually are~\footnote{However, there are examples of negatively associated discrete measures including determinantal ones  (see \cite[Theorem 6.5]{Lyons03}).}.

In order to unify the approach to matter in hand 
and provide more examples of pp comparable to Poisson pp, 
we define two more classes of pp:
{\em weakly sub-Poisson} --- as pp having  
both void probabilities and factorial moment
measures smaller than the Poisson pp with same mean measure,
and {\em weakly super-Poisson} --- as having these characteristics
larger than the Poisson pp with same mean measure.
It is almost straightforward to see that this new classification is indeed weaker than
sub- and super-poissonianity based on the $dcx$ ordering. 
We prove that it is also weaker than association:
positive association implies weak
super-poissonianity, while negative association implies 
weak sub-poissonianity. A good news is that 
permanental and determinantal pp can be proved to be weakly
super- and sub-Poisson respectively. Also, as it turns out many of the 
results can be proven under these weaker assumptions of weakly sub-Poisson or super-Poisson 
than association or $dcx$ ordering. 

\paragraph{Paper organization} The necessary notions, notations
and basic facts are introduced and recalled in Section~\ref{s.Notation}. In
Section~\ref{s.sub-super-Poisson}, we define classes of strongly
and weakly sub- and super-Poisson pp and, as a main result, 
we prove that weak sub- or super-poissonianity is implied by negative or positive
association, respectively. We study the perturbed-lattice pp in Section~\ref{sec:ex} 
and determinantal and permanental pp in Section~\ref{sec:det.perm}.
In Section~\ref{s.Applications}, we discuss some further theoretical implications (especially percolation) of 
the presented ideas as well as their connections to other stochastic geometric models and the modelling applications. 
Lemma~\ref{lem:MeesterSh-like}, which is of independent interest and 
 used in this paper for showing $dcx$ ordering of perturbed lattices and determinantal and permanental point processes
(on mutually disjoint  simultaneously observable sets) is  proved in the Appendix.

\section{Notions, notation and basic facts}
\label{s.Notation}
\paragraph{Point processes}
We assume the usual framework for random measures and point processes
on $d$-dimensional Euclidean space $\mR^d$ ($d\ge1$), 
where these are considered as random elements 
on  the  space  $\mathbb{M}(\mR^d)$
of non-negative Radon
measures on~$\mR^d$ (cf \cite{Kallenberg83}).
A point process (pp) $\Phi$ is
simple if a.s. $\Phi(\{x\}) \leq 1$ for all $x \in \mR^d$. 
We denote by $\nu(B)=\pr{\Phi\cap B=\emptyset}$,
the void probabilities of pp  $\Phi$ and
by $\alpha^{(k)}(\cdot)$, the {\em factorial moment measure} of~$\Phi$.
Recall that for simple pp, $\alpha^{(k)}(B_1\times\ldots\times B_k)
=\sE(\prod_{i=1}^k \Phi(B_i))$ for pairwise disjoint bounded Borel subsets(bBs) $B_i$  ($i=1,\ldots, k$).
The {\em $k\,$th  joint intensity}, $\rho^{(k)}:(\mR^{d})^{k} \to
[0,\infty)$ is the density (if it exists) of $\alpha^{(k)}(\cdot)$
with respect to the Lebesgue measure $\md x_1\dots\md x_k$.
Recall that the  joint intensities $\rho^{(k)}$, $k\ge 1$
characterize the distribution of a pp.
The above facts remain true even when the densities $\rho^{(k)}$
are considered with respect to  $\prod_{i=1}^k\mu(\md x_k)$ for an
arbitrary Radon measure $\mu$ on $\mR^d$.  
As always, a pp or a random measure on $\mR^d$ is said to be {\em stationary} if its distribution is invariant with respect to translation by
vectors in~$\mR^d$.

\paragraph{Directionally convex ordering}
A Lebesgue-measurable function $f:\mR^k \rar \mR$ is said to be {\em
  directionally convex}~($dcx$) if for every $x \in \mR^k,
\epsilon,\delta > 0, i,j \in \{1,\ldots,k\}$, we have that
$\Delta_{\epsilon}^i \Delta_{\delta}^jf(x) \geq 0$, where
$\Delta_{\epsilon}^if(x) := f(x + \epsilon e_i) - f(x)$
is the discrete differential
operator, with  $\{e_i\}_{1 \leq i \leq k}$
denoting  the canonical basis vectors for $\mR^k$. 
 We abbreviate {\em  increasing} and $dcx$ by $idcx$ and {\em
   decreasing} and $dcx$ by $ddcx$ (see \cite[Chapter 3]{Muller02}). 
For real-valued random vectors of the
same dimension $X$ and $Y$, {\em $X$ is said to be less than $Y$ in
$dcx$ order} (denoted  $X \leq_{dcx} Y$)
if $\sE(f(X)) \leq \sE(f(Y))$ for all $f$ $dcx$
such that both the expectations are finite. 
For two pp on $\mR^d$, one says that  $\Phi_1(\cdot) \leq_{dcx} \Phi_2(\cdot)$, if for any $B_1, \ldots, B_k$ bBs in $\mR^k$,
$(\Phi_1(B_1),\ldots,\Phi_1(B_k)) \leq_{dcx}
(\Phi_2(B_1),\ldots,\Phi_2(B_k))$;
cf~\cite{snorder}.
The definition is similar for other orders, i.e., those defined by $idcx,\allowbreak ddcx$ functions.
It is enough to verify the above
conditions for $B_i$ mutually disjoint.
In order to avoid technical difficulties, we will
consider here only pp whose  {\em mean measures} $\sE(\Phi(\cdot))$
are Radon (finite on bounded sets). For such pp, $dcx$ order is a
transitive order. Due to the fact that each $dcx$ function
can be monotonically approximated by $dcx$ functions $f_i(\cdot)$ which satisfy $f_i(x)=O(||x||_\infty)$ at infinity, where $||x||_\infty$ is the $L_\infty$ norm on the Euclidean space; cf.~\cite[Theorem~3.12.7]{Muller02}.

It is easy to see that 
$\Phi_1(\cdot) \leq_{dcx} \Phi_2(\cdot)$ implies the {\em equality
of their  mean measures}: $\sE(\Phi_1(\cdot)) =\sE(\Phi_2(\cdot))$.
Moreover, as shown in~\cite{snorder}, 
higher-order moment measures are non-decreasing in $dcx$ order on pp, provided
they are $\sigma$-finite%
\footnote{$\sigma$-finiteness condition is missing in~\cite{snorder}; see~\cite[Prop.~4.2.4]{Yogesh_thesis} for the
 correction}. In addition,  $dcx$ ordering allows  to compare
also the void probabilities as stated in the following  
new result:
\begin{prop}
\label{Prop:voids_pp}
Denote by $\nu_1(\cdot),\nu_2(\cdot)$ the void probabilities of pp
$\Phi_1$ and $\Phi_2$ on $\mR^d$ respectively.
If  $\Phi_1 \leq_{ddcx} \Phi_2$ then
$\nu_1(B) \le \nu_2(B)$
for all bBs $B\subset\mR^d$.
\end{prop}
\begin{proof}
This follows  directly from the definition of $dcx$
ordering of pp, expressing  $\nu_j(B)=\sE(f(\Phi_j(B)))$, $j=1,2$, with the
function $f(x)=\max(0,1-x)$ that is decreasing and convex (so
$ddcx$ in one dimension). \qed
\end{proof}
In particular, the latter result implies
ordering of all contact distribution functions (empty space functions)
for pp comparable in $dcx$ order and not having 
fixed atoms~\footnote{Satisfying $\Pr\{x\not\in\Phi\}=1$ for all
  $x\in\mR^d$.}.
We see in the joint comparison of moment measures and void
  probabilities of pp having equal mean measures,  
a new tool for comparison of their clustering properties, weaker than $dcx$ order but more easy to verify.

\begin{sloppypar} 
\paragraph{Positive and negative association}
Denote by  $\COV{XY}=\EXP{XY}-\EXP{X}\EXP{Y}$ covariance of
random variables $X,Y$.
A point process  $\Phi$ 
is called  {\em associated} if 
$\COV{f(\Phi(B_1),\ldots,\Phi(B_k)),g(\Phi(B_1),\ldots,\Phi(B_k))}\allowbreak \geq 0$
for any finite collection of bBs $B_1,\ldots,B_k\subset\mR^d$ and
$f,g$ continuous and increasing functions taking values in $[0,1]$;
cf~\cite{BurtonWaymire1985}. This property is also called {\em
  positive association}, or the {\em FKG property}.
The theory for the opposite property is more tricky, cf~\cite{pemantle00},
but one can call 
$\Phi$ {\em negatively associated} if
$\COV{f(\Phi(B_1),\ldots,\Phi(B_k)),\allowbreak g(\Phi(B_{k+1}),\ldots,\Phi(B_l))}\allowbreak \le 0$
for any finite collection of bBs $B_1,\ldots,B_l\subset\mR^d$ 
such that $(B_1\cup\dots\cup  B_k)\cap(B_{k+1}\cup\dots \cup
B_{l})=\emptyset$ and
$f,g$ increasing functions;
Both definitions can be straightforwardly extended to
 random measures.
\end{sloppypar}

\section{Comparison of clustering to Poisson pp}
\label{s.sub-super-Poisson}
We  call a pp
{\em sub-Poisson} (respectively {\em
  super-Poisson}) if it is smaller (larger) in $dcx$ order than the
Poisson pp (necessarily of the same mean measure). 
(More precisely we should have called these processes  $dcx$-sub-Poisson or
$dcx$-super-Poisson pp, but we omit the word $dcx$ for simplicity.)
Examples of such pp
are given in Section~\ref{sec:ex}. 
A weaker notion of sub-  and
super-poissonianity  can be defined when comparing only moment measures or void
probabilities. 
Bearing in mind that Poisson pp
can be characterized as having void probabilities of the form
$\nu(B)=\exp(-\alpha(B))$, where $\alpha(\cdot)$ is its mean measure,
we say that a pp $\Phi$
is {\em weakly sub-Poisson in the sense of void probabilities}
({$\nu$-weakly sub-Poisson}) if
\begin{equation}\label{e.nu-weakly-sub-Poisson} 
\pr{\Phi(B)=0}\le
e^{-\EXP{\Phi(B)}}
\end{equation}
for all Borel sets $B \subset \mR^d$. Similarly, we say that a pp
$\Phi$ is {\em weakly sub-Poisson in the sense of moment measures}
({$\alpha$-weakly sub-Poisson}) if \begin{equation}\label{e.alpha-weakly-sub-Poisson}
\alpha^{(k)}(B_i) \le \prod_{i=1}^k\alpha^{(1)}(B_i) = \prod_{i=1}^k\EXP{B_i} \,
\end{equation} for all mutually disjoint bBs $B_i\subset\mR^d$.  When the inequalities
in~(\ref{e.nu-weakly-sub-Poisson}) and (\ref{e.alpha-weakly-sub-Poisson}) are reversed, we will say that $\Phi$
is $\nu$-weakly super-Poisson or $\alpha$-weakly super-Poisson respectively.

Finally, we will say that $\Phi$ is {\em weakly sub-Poisson}
if $\Phi$ is $\alpha$-weakly sub-Poisson and $\nu$-weakly sub-Poisson.
Similarly, we define {\em weakly super-Poisson} pp.
Examples of weakly sub- and super-Poisson pp are given in
Section~\ref{sec:det.perm}.

The fact that $dcx$ ordering implies ordering of moment
 measures and void probabilities
lend credence to our usage of the terms weak
sub- and super-Poissonianity. 
Interestingly, these inequalities  are also implied by 
negative and positive association. The following result
is a key observation in this matter.

\begin{prop}\label{p.association}
Consider point process $\Phi$ with 
Radon mean measure $\alpha(\cdot)=\EXP{\Phi(\cdot)}$.
If $\Phi$ is simple, has Radon  second-order factorial
moment measure $\alpha^{(2)}(\cdot)$ and
\begin{equation}\label{eq:2_void-inequality}
\Pr\{\,\Phi(B_1)=0, \Phi(B_2)=0\,\}\le
\Pr\{\,\Phi(B_1) = 0\,\}\Pr\{\,\Phi(B_2)=0\,\},
\end{equation}
for any two disjoint bBs $B_1$ and $B_2$,
then  $\Phi$ is $\nu$-weakly sub-Poisson.

If  the mean measure $\alpha(\cdot)$ of $\Phi$  
is diffuse (without atoms)
and $\Phi$ satisfies~(\ref{eq:2_void-inequality}) 
with the reversed inequality ($\ge$)
for any two disjoint bBs $B_1$ and $B_2$,
then $\Phi$ is $\nu$-weakly super-Poisson.
\end{prop}

\begin{proof}
\begin{sloppypar}

Define a set function $Q(B)=-\log(\Pr\{\Phi(B)=0\,\})$.
Regarding the  first statement,
it is immediate to see that $Q$ is non-negative and, 
under assumption~(\ref{eq:2_void-inequality})
super-additive; i.e., 
for any finite
$k \ge1$ and any pairwise disjoint bBs~$B_j$,
$j=1,\ldots,k$  
$Q(B_1\cup\ldots\cup B_k)\ge \sum_{j=1}^k Q\left( B_j\right)$. 
In order to prove the result, we need to  show that $Q(B)\ge \alpha(B)$, for
any bBs~$B$. 
To this regard note by  the super-additivity of $Q$ that for any
bBs~$B$
\end{sloppypar}
\begin{equation}\label{eq:tildeQsub}
Q(B)=\sup_{J}\sum_{j \in J}Q(B_j)\,,
\end{equation}
where the ``sup'' is taken over all {\em finite} partitions of $B$ into
bBs $B_i$. Moreover, for any bBs~$B$ 
\begin{eqnarray*}\nonumber
\Pr\{\,\Phi(B)=0\,\}&=&1-\EXP{\Phi(B)}+\EXP{\Phi(B)\1(\Phi(B)\ge2)} - \Pr\{ \Phi(B) \ge 2\} \nonumber \\
&\le& 1-\EXP{\Phi(B)}+\EXP{\Phi(B)(\Phi(B)-1)^+}\nonumber\\
&=& 1-\alpha(B)+\alpha^{(2)}(B\times B)\,
\end{eqnarray*}
and hence 
$Q(B)=-\log(\Pr\{\,\Phi(B)=0\,\})\ge \alpha(B)-\alpha^{(2)}(B\times B)$.
Consequently, by~(\ref{eq:tildeQsub}), for any bBs $B$
$$Q(B)\ge \sup_{J}\sum_{j \in J}\left(\alpha(B_j)-\alpha^{(2)}(B_j\times
B_j)\right)=\alpha(B)-\inf_{J}\sum_{j\in J}\alpha^{(2)}(B_j\times
B_j)\,,$$
due to finiteness of all terms.
In order to complete the proof it is enough to show that the 
``inf'' term is equal to zero. To this regard, for a given $\epsilon>0$
define $\Delta_B^\epsilon=\{B\times B\owns(x,y):  |(x,y)-(z,z)|\le
\epsilon\; {\rm for\; some}\; z\in B\}$.
Note that $\Delta_B^\epsilon $  can be seen as some neighborhood of the
intersection of the diagonal with $B\times B$.
Note also that for any $\epsilon>0$ there exits
a suitable fine partition $I$ of $B$ 
such that 
$\sum_{j\in J}\alpha^{(2)}(B_j\times B_j)\le
\alpha^{(2)}(\Delta_B^\epsilon)$.
(For example, take a finite  coverage of $B$ by balls of radius
$\epsilon$, which exists by local-compactness of the space,
and refine it to have disjoint partition of $B$.) 
By the local finiteness and $\sigma$-additivity of  $\alpha^{(2)}$,
$\lim_{\epsilon\to0}
\alpha^{(2)}(\Delta_B^\epsilon)=
\alpha^{(2)}(\{(z,z):z\in B\})=0$,
where the last equality follows from the assumption that $\Phi$
is simple. This completes the proof of the first statement.

For the second statement, we will show that $Q(B) \le \alpha(B)$.
To this regard, note that the reversed inequality
in~(\ref{eq:2_void-inequality}) implies that  $Q(\cdot)$ is
sub-additive and consequently, for any bBs $B$,
\begin{equation}\label{eq:tildeQsuper}
Q(B)=\inf_{J}\sum_{j \in J}Q(B_j)\,,
\end{equation}
where ``inf'' is over all finite partitions of $B$.
Moreover, observe that 
$\Pr\{\,\Phi(B)=0\,\}\ge 1-\alpha(B)$
and that, for $0\le x\le\epsilon$, 
$-\log(1-x) \le x(1+\delta(\epsilon))$, where 
$\delta(\epsilon)=\epsilon/(2(1-\epsilon)^2$,
which can be shown by the Taylor expansion with Lagrange form of the
remainder term of order~2. Since $\alpha(\cdot)$ is diffuse, 
for any $\epsilon>0$ there exists 
a partition $J$ of bBs~$B$ 
such that $\alpha(B_j)\le\epsilon$ for all $j\in J$. 
For such a partition $J$,
$$Q(B)\le \sum_{j\in J}-\log(1-\alpha(B_j))\le
\alpha(B)(1+\delta(\epsilon))\,.$$
The proof follows from the observation that $\delta(\epsilon)\to 0$
when $\epsilon\to 0$. \qed
\end{proof}

\begin{cor}\label{c.napa_vs_Poisson}
A negatively associated, simple pp with a Radon mean measure 
is weakly sub-Poisson.
A (positively) associated pp with a Radon, diffuse  mean measure 
is weakly super-Poisson.
\end{cor}
\begin{proof}
Inequality~(\ref{e.alpha-weakly-sub-Poisson}) or its inverse (i.e.;
$\alpha$-weak sub- or super-poissonianity) follows directly from
negative association or association, respectively.
The $\nu$-weak sub- or super-poissonianity follows from
Proposition~\ref{p.association}. Indeed, inequality
(\ref{eq:2_void-inequality}) or its inverse can be
derived easily from negative association or association, respectively.
Moreover, note by~(\ref{e.alpha-weakly-sub-Poisson}),
that any factorial  moment measure $\alpha^{(n)}(\cdot)$
of a simple,  $\alpha$-weakly sub-Poisson pp with Radon
mean measure is also Radon. This completes the proof. \qed
\end{proof}

In fact, sub-Poissonianity (or negative association provided  the
aforementioned regularity of pp)  implies
something stronger than $\alpha$-weak sub-Poissonianity.  
Namely, we have that, $\alpha^{(k+l)}(\cdot) \leq \alpha^{(k)}(\cdot)\alpha^{(l)}(\cdot)$ for integers $k,l \geq 0$. 
Similarly  super-Poissonianity (or positive association provided the
aforementioned regularity of pp)
implies the reverse inequality. Further justification for negative association as a measure of sparsity will be seen in \cite{Adler12} where
it is shown that Palm measure of a negatively associated pp is
``stochastically weaker'' than that of the original pp. In particular,
the void probability 
increases for the Palm measure.

\subsubsection*{A counterexample.}
Let us finally remark existence of negatively associated pp which are
not sub-Poisson (neither in $dcx$ nor weakly). Our
counterexample is not a simple pp,  which shows also that  this
latter assumption cannot be relaxed in Corollary~\ref{c.napa_vs_Poisson}.
In this regard,
for a given fixed
integer $k$ consider a discrete subset  $\{x_1,\ldots,x_k\}$ of the space
and a point process $\Phi$ supported on this set, such that the
vector $(N_1,\ldots,N_k)$, with $N_i=\Phi(\{x_i\})$
has the {\em permutation distribution} of the vector 
$(0,1,\ldots,k-1)$,  i.e., it takes as values all $k!$ permutations
of this vector with equal probabilities, each being $1/k!$.
By~\cite[Theorem~2]{JoagDev1983}  $(N_1,\ldots,N_k)$ and hence $\Phi$ 
is negatively associated. Note that $N_i$ is uniform random variable
on $\{0,1,\ldots,k-1\}$.  Thus it has mean $\EXP{N_i}=(k-1)/2$,
void probability $\Pr\{\,N_i=0\,\}=1/k$ and variance $(k^2-1)/12$.
Note that for sufficiently large $k$ we have $1/k>e^{-(k-1)/2}$
and $(k^2-1)/12>(k-1)/2$; i.e., the void probability and the variance
of  $N_i$ are  larger than these of 
Poisson variable of mean $(k-1)/2$. Consequently $\Phi$ is not
sub-Poisson (in $dcx$ sense) and not  $\nu$-weakly sub-Poisson.

\section{Perturbed lattices and point processes}
\label{sec:ex}
\hfill\newline
It was observed in~\cite{snorder} that Poisson-Poisson cluster pp,
L\'{e}vy  based Cox pp, Ising-Poisson cluster pp are super-Poisson pp.
In this section, we present more examples of pp, which are $dcx$
comparable to Poisson pp.
We begin with  a
general model of a perturbation of a pp and prove 
our key result on the $dcx$ ordering of such pp.

\subsection{Perturbation operator}
\label{ssec:pl}

Let $\Phi$ be a pp on $\mR^d$ and $\cN(\cdot,\cdot)$,
$\cX(\cdot,\cdot)$ be two probability kernels from $\mR^d$ to
non-negative integers $\mZ^+$ and $\mR^d$, respectively.
Consider the following {\em independently marked}
version of the pp $\Phi$, $\tilde\Phi^{pert}=\{(X,N_X,{\mathbf Y}_X) \}_{X\in\Phi}$
where given $\Phi$:\vspace{-2ex}
\begin{itemize}
\item $N_X$,  $X\in\Phi$ are independent, non-negative
integer-valued random variables with distribution
$\pr{N_X\in\cdot\,|\,\Phi}=\cN(X,\cdot)$,
\item  ${\mathbf Y}_{X}=(Y_{iX} : i =1,2,\ldots)$, $X\in\Phi$ are
  independent vectors of  i.i.d.  elements of  $\mR^d$, with
  $Y_{iX}$'s having the conditional distribution $\pr{Y_{iX}\in\cdot\,|\,\Phi}=\cX(X,\cdot)$,
\item the random elements
$N_X,{\mathbf Y_X}$ are independent for all $X\in\Phi$.
\end{itemize}
Consider the following subset of $\mR^d$
\begin{equation}
\label{defn:perturbed_pp}
\Phi^{pert}= \bigcup_{X \in \Phi} \bigcup_{i=1}^{N_X} \{ X + Y_{iX}\}\,,
\end{equation}
where the inner sum is interpreted as $\emptyset$ when $N_X=0$.
The set $\Phi^{pert}$ can  (and will) be considered as a pp
on $\mR^d$ provided it is locally finite. In what follows, in
accordance with our general assumption for this article, we will assume that the mean measure of $\Phi^{pert}$ is locally finite (Radon measure)
\begin{equation}
\label{e:perturbation-Radon}
\int_{\mR^d} n(x)\cX(x,{B-x})\,  \al(dx)< \infty, \quad
\text{for all bBs~$B\subset\mR^d$},
\end{equation}
where  $\al(\cdot)$ is the mean measure of the pp $\Phi$ and
$n(x)=\sum_{k=1}^\infty 
k\cN(x,\{k\})$ is the mean value of the distribution $\cN(x,\cdot)$.

The pp $\Phi^{pert}$ can be seen as
{\em independently replicating and
translating points from the pp $\Phi$}, with the  number of
replications of the point  $X\in\Phi$ having distribution
$\cN(X,\cdot)$ and the independent  translations of these replicas
from $X$  by vectors having distribution $\cX(X,\cdot)$. For this
reason, we call $\Phi^{pert}$ a {\em perturbation} of~$\Phi$
driven by  the {\em replication kernel} $\cN$ and the {\em translation kernel}
$\cX$.

An important observation for us is that the {\em operation of perturbation of $\Phi$ is $dcx$ monotone with respect to the replication kernel} in the following sense.

\begin{prop}\label{p.pert-lattice}
Consider a pp $\Phi$ with Radon mean measure $\al(\cdot)$
and its two perturbations
$\Phi_j^{pert}$ $j=1,2$ satisfying
condition~(\ref{e:perturbation-Radon}), having the same
translation kernel $\cX$ and possibly different replication kernels
$\cN_j$, $j=1,2$, respectively.
If $\cN_1(x,\cdot)\leq_{cx}\cN_2(x,\cdot)$ (convex
ordering of the conditional distributions of the number of replicas)
for $\al$-almost all $x\in\mR^d,$ then
$\Phi^{pert}_1\leq _{dcx} \Phi^{pert}_2$.
\end{prop}
\begin{proof}
We will consider some particular coupling of the two perturbations
$\Phi^{pert}_j$,  $j=1,2$.
Given  $\Phi$ and ${\mathbf Y}_X=(Y_{iX}:i=1,\ldots)$ for
each  $X\in\Phi$, let $\Phi_{jX} = \bigcup_{i=1}^{N^j_{X}} \{ X + Y_{iX} \}$,
where $N^j_X$ has distribution $\cN_j(X,\cdot)$, $j=1,2$,
respectively.
Thus  $\Phi^{pert}_j=\sum_{X\in\Phi}\Phi_{jX},$ $j=1,2$
are the two considered perturbations.
Note that given $\Phi$, $\Phi^{pert}_j$ can be seen as
independent superpositions of $\Phi_{jX}$ for $X \in \Phi.$ Hence,
by \cite[Proposition 3.2(4)]{snorder} (superposition preserves $dcx$
order) and \cite[Theorem 3.12.8]{Muller02} (weak and $L_1$ convergence
jointly preserve $dcx$ order),
it is enough to show that conditioned on $\Phi$, $\Phi_{1X} \leq_{dcx}
\Phi_{2X}$ for every $X \in\Phi$.
In this regard, given $\Phi$, consider $X\in\Phi$ and let
$B_1,\ldots,B_k$ be mutually disjoint bBs and $f : \mR^k
\to \mR$, a $dcx$ function. Define a real valued function
$g:\mZ \to \mR$, as
$$g(n):= \EXP{f\Bigl(\text{sgn}(n)\sum_{i=1}^{|n|}(\1[Y_{iX} \in B_1 - X],\ldots,
\1[Y_{iX} \in B_k - X])\Bigr) \bigg| \Phi}\,,$$
where $\text{sgn}(n)= \frac{n}{|n|}$ for $n \neq 0$ and $\text{sgn}(0) = 0.$ By Lemma~\ref{lem:MeesterSh-like}, $g(\cdot)$ is a convex function
on~$\mZ$ and by Lemma~\ref{lem:discrete-convex}
it can be extended to a convex function $\tilde g(\cdot)$ on~$\mR$.
Moreover,
$\EXP{\tilde g(N^j_{X}) | \Phi} =\EXP{g(N^j_{X}) | \Phi}\allowbreak =
\sE\Bigl(f(\Phi_{jX}(B_1),\ldots,\Phi_{jX}(B_k)) | \Phi\Bigr)$
for $j=1,2$.
Thus, the result follows from the assumption $N_{X}^1\le_{cx}N_{X}^2$. \qed
\end{proof}

\begin{rem}
The above proof remains valid for an extension
of the perturbation model in which
the distribution $\cX(X,\cdot)$ of the translations $Y_{iX}$
depends not only on the location of the point $X\in\Phi$ but also on the
entire configuration $\Phi$; $\cX(X,\cdot)=\cX(X,\Phi, \cdot)$,
provided condition~(\ref{e:perturbation-Radon}) is replaced by finiteness of
$\int_{\mM^d}\int_{\mR^d}n(x)\cX(x,\phi,B-x)\,C(d(x,\phi))$,
where $C(d(x,\phi))$ is the Campbell measure of $\Phi$.
\end{rem}

\subsection{Examples}
\label{ss.pert-Examples}
\paragraph{Perturbed Poisson pp}
 Let $\Phi$ be a (possibly
  inhomogeneous) Poisson pp of mean measure $\alpha(dx)$ on $\mR^d$.
Let
$\cN(x,\cdot)=\varepsilon_1={\1(1\in\cdot)}$ be the Dirac measure
on $\mZ^+$ concentrated at~1 for all $x\in\mR^d$ and assume
 an  arbitrary translation kernel $\cX$ satisfying
$\alpha^{pert}(A)=\int_{\mR^d}\cX(x,\allowbreak
A-x)\,\alpha(dx)<\infty$ for all
bBs $A$. Then  by the displacement theorem for Poisson pp,
$\Phi^{pert}$ is also a Poisson pp with mean measure $\alpha^{pert}(dx)$.
Assume {\em any} replication kernel $\cN_2(x,\cdot)$,
with mean number of replications $n_2(x)=\allowbreak\sum_{k=1}^\infty
k\cN_2(x,\{k\})=1$ for all $x\in\mR^d$. Then, by the Jensen's inequality and
Proposition~\ref{p.pert-lattice}, one obtains a super-Poisson pp
$\Phi_2^{pert}$. In the special case, when $\cN_2(x,\cdot)$ is
the Poisson distribution with  mean~1 for all $x\in\mR^d$,
$\Phi_2^{pert}$ is a Poisson-Poisson cluster pp which is
a special case of a Cox (doubly stochastic Poisson) pp with
(random) intensity measure  $\Lambda(A) = \sum_{X \in
  \Phi}\cX(x,A-x)$. The fact that it is super-Poisson was already observed
in~\cite{snorder}. Note that for a general distribution of  $\Phi$,
its perturbation $\Phi_2^{pert}$ is also a Cox pp of the intensity
$\Lambda$ given above.  
\paragraph{ Perturbed lattice pp}
Assuming a deterministic lattice $\Phi$ (e.g. $\Phi=\mZ^d$)
gives rise to the perturbed lattice pp of the type considered
in~\cite{Sodin04}. Surprisingly enough, starting
from such a $\Phi$, one can also construct a Poisson pp and
both super- and sub-Poisson perturbed pp. In this regard, assume for simplicity that $\Phi=\mZ^d$, and the translation kernel
$\cX(x,\cdot)$ is uniform on the unit cube $[0,1)^d$.
Let $\cN(x,\cdot)$ be the Poisson distribution with mean~$\lambda$ ($Poi(\lambda)$). It is easy to see that such a perturbation $\Phi^{pert}$ of the lattice $\mZ^d$ gives rise to a homogeneous Poisson pp with intensity~$\lambda$.
\subsubsection{Sub-Poisson perturbed lattices.}
\label{sec:sub_poisson_lattice}
Assuming for  $\cN_1$ some distribution convexly ($cx$)  smaller than $Poi(\lambda),$
one obtains a sub-Poisson perturbed lattice pp.
Examples are {\em hyper-geometric} $H\!Geo(n,m,k)$,
$m,k\le n$, $km/n=\lambda$
and {\em binomial}  $Bin(n,\lambda/n)$, $\lambda\le n$ distributions%
~\footnote{$Bin(n,p)$ has probability mass function
  $p_{Bin(n,p)}(i)={n\choose i}p^i(1-p)^{n-i}$ ($i=0,\ldots,n$).
$H\!Geo(n,m,k)$ has probability
  mass function $p_{H\!Geo(n,m,k)}(i)={m\choose i}{n-m\choose
    k-i}/{n\choose k}$ ($\max(k-n+m,0)\le i\le m$).},
which can be ordered as follows:
\begin{equation}
\label{eqn:sub_Poisson_rvs}
H\!Geo(n,m,\lambda n/m)\le_{cx} Bin(m,\lambda/m)\le_{cx}
Bin(r,\lambda/r)\le_{cx} Poi(\lambda),
\end{equation}
for $\lambda\le m\le \min(r,n)$;
cf.~\cite{Whitt1985}%
\footnote{One shows the
logarithmic concavity of the ratio of the respective probability mass
functions, which  implies increasing convex order and, consequently,
$cx$ provided the distributions have the same  means.}.
Specifically, taking  $\cN_1(x,\cdot)$ to be Binomial
$Bin(n,\lambda/n)$ for
$n\ge\lambda$, one obtains
a $dcx$ monotone increasing family of sub-Poisson pp.
Taking $\lambda=n=1$ (equivalent to $\cN(x,\cdot)=\varepsilon_1$),
one obtains a {\em simple perturbed lattice} that is $dcx$ smaller
than the Poisson pp of intensity~1.

\subsubsection{Super-Poisson perturbed lattices.}
\label{sec:sup_poisson_lattice}
Assuming for  $\cN_2$ some distribution convexly larger than $Poi(\lambda),$
one obtains a super-Poisson perturbed lattice. Examples are
{\em negative binomial} ${N\!Bin}(r,p)$ distribution
with $rp/(1-p)=\lambda$ and {\em geometric} $Geo(p)$ distribution with
$1/p-1=\lambda$~%
\footnote{$p_{Geo(p)}(i)=p(1-p)^{i}$,  $p_{N\!Bin(r,p)}(i)=
{r+i-1\choose i}p^i(1-p)^r$.},
which can be ordered in the following way:
\begin{eqnarray}
Poi(\lambda) & \le_{cx} & N\!Bin(r_2,\lambda/(r_2+\lambda))
\le_{cx} N\!Bin(r_1,\lambda/(r_1+\lambda)) \no \\
\label{eqn:sup_Poisson_rvs} & \le_{cx} & Geo(1/(1+\lambda)) 
\le_{cx}\sum_{j}\lambda_j\;Geo(p_j) 
\end{eqnarray}
with $r_1\le r_2$,
$0\le \lambda_j\le 1$, $\sum_j\lambda_j=1$ and
$\sum_j\lambda_j/p_j=\lambda+1$,
where the largest distribution above is a mixture
of geometric distributions having mean~$\lambda$;
cf.~\cite{Whitt1985}.
Specifically, taking  $\cN_2(x,\cdot)$ to be negative binomial
$N\!Bin(n,\lambda/(n+\lambda))$ for
$n=1,\ldots$ one obtains
a $dcx$ monotone decreasing  family of super-Poisson pp.
Recall that $N\!Bin(r,p)$
is  a mixture of $Poi(x)$  with parameter $x$ distributed
as a gamma distribution with scale parameter $p/(1-p)$
and shape parameter~$r$.

From~\cite[Lemma~2.18]{MeesterSh93}, we know that any {\em  mixture of Poisson distributions} having mean $\lambda$ is $cx$ larger than $Poi(\lambda)$.
Thus, the super-Poisson perturbed lattice with such a replication kernel
(translation kernel being the uniform distribution) 
again gives rise to a Cox pp.

\subsubsection{Associated point processes:}
\label{sec:assoc_pp}

From~\cite[Th.~5.2]{BurtonWaymire1985}, we know that any {\em Poisson
center cluster pp} is (positively) associated. This is a generalization of our
perturbation~(\ref{defn:perturbed_pp}) of a Poisson pp $\Phi$
(cf. Section~\ref{ss.pert-Examples}) having form
$\Phi^{cluster}=\sum_{X\in\Phi}\{X+\Phi_X\}$ with $\Phi_X$ being
arbitrary i.i.d. (cluster) point measures.  Other examples of
associated pp given in~\cite{BurtonWaymire1985} 
are Cox pp with intensity measures being associated.
(It is easy to see by Jensen's inequality that all Cox pp are $\nu$-weakly super-Poisson.)

It is easy to see that the pp formed by throwing $n$ i.i.d. points in a bounded region forms
a negatively associated pp. Further, one can show that independent
superposition of negatively associated pp is a negatively associated pp. Hence, simple perturbed lattices
(cf. Section \ref{sec:sub_poisson_lattice}) are negatively associated.

\section{Determinantal and permanental point processes} 
\label{sec:det.perm}

In this section, we focus on spatial determinantal and
permanental pp. We will show that they are, respectively, weakly sub-
and super-Poisson pp. Some partial $dcx$ comparison of these pp
with respect to Poisson pp, namely on mutually disjoint, simultaneously observable sets, will be proved as well.

\subsection{Definition}
\label{sss.Int_ker}
To make the paper more self-contained,
we will recall a general framework from \cite[Chapter~4]{Ben09},
which allows  us to study ordering of determinantal and permanental pp
more explicitly; see also \cite{Ben06} for a quick
introduction to these pp.

Let $K:\mR^{d}\times\mR^d \to \mathbb{C}$ (where $\mathbb{C}$ are complex numbers) be a {\em locally square-integrable} kernel, with respect to  $\mu^{\otimes 2}$ on $\mR^{2d}$~\footnote{i.e., $\int_D\int_{D}|K(x,y)|^2\,\mu(dx)\mu(dy)<\infty$ for every compact $D\subset\mR^d$}. Then $K$ defines an associated integral operator $\mathcal{K}_D$ on $L^2(D,\mu)$ as $\mathcal{K}_Df(x)=\int_DK(x,y)f(y)\,\mu(dy)$ for complex-valued, square-integrable $f$ on $D$ ($f\in L^2(D,\mu)$).
This operator is compact and hence its  spectrum is discrete.
The only possible accumulation point is $0$ and every non-zero
eigenvalue has finite multiplicity.
Assume moreover that for each compact $D$ the operator $\mathcal{K}_D$
is {\em Hermitian}~\footnote{i.e., $\int_D\overline
{f(x)}\mathcal{K}_Dg(x)\,\mu(dx)=\int_D\overline
{g(x)}\mathcal{K}_Df(x)\,\mu(dx)$ for all $f,g\in L^2(D,\mu)$},
{\em positive semi-definite}~\footnote{i.e.,
 $\int_D\overline {f(x)}\mathcal{K}_Df(x)\,\mu(dx)\ge 0$}, and
{\em trace-class}; i.e.,
$\sum_{j}|\lambda_j^D|<\infty$, where $\lambda_j^D$ denote the
  eigenvalues of $\mathcal{K}_D$. By the positive semi-definiteness of
  $\mathcal{K}_D,$ these eigenvalues are non-negative.
\remove{
Further, one can show (cf.~\cite[Lemma~4.2.2]{Ben09}) that for each compact $D$, there exists a ``version''  $K_D(x,y)$
of the kernel $K$, defined on $D'\subset D$ such that $\mu(D\setminus
D')=0$, having the  same associated operator $\mathcal{K}_D$ on $L^2(D,\mu)$
~\footnote{i.e. $K_D(x,y)=K(x,y)$ for $\mu^{\otimes 2}$ almost all
$x,y\in D$}, which
is Hermitian and positive semi-definite.\footnote{%
Recall, a kernel $K(x,y)$ is Hermitian if
$K(x,y)=\overline {K(y,x)}$ for all $x,y\in\mR^d$, where $\overline z$ is the complex
conjugate of $z\in\mathbb{C}$. It is positive semi-definite
$\sum_{i=1}^k\overline z_i \sum_{j=1}^k K(x_i,x_j) z_j\ge 0$
for all $z_i\in\mathbb{C}$, $i=1,\ldots,k$, $k \geq 1.$}
Specifically, one can take $K_D(x,y) = \sum_j \lam_j^D
\phi_j^D(x) \overline{\phi_j^D(y)}$ where $\phi_j^D(\cdot)$ are the
corresponding normalized eigenfunctions of $\mathcal{K}_D$.
}

\paragraph{Determinantal pp}
A simple pp on $\mR^d$
is said to be a {\em determinantal pp} with a kernel $K(x,y)$ with respect to
a Radon measure $\mu$ on $\mR^d$ if the  joint intensities
of the pp with respect to the product measure $\mu^{\otimes k}$ satisfy
$\rho^{(k)}(x_1,\ldots,x_k) = \det \big( K(x_i,x_j) \big)_{1
\leq i,j \leq k}$ for
all $k$, where $\big(a_{ij}\big)_{1 \leq i,j \leq k}$ stands for a
matrix with entries $a_{ij}$ and $\det \big( \cdot \big)$ denotes the
determinant of the matrix. Note that the mean measure of the
determinantal pp (if it exists) is  equal to
$\alpha(\cdot)=\int_{\cdot}K(x,x)\,\mu(dx)$.
Assuming that the kernel $K$ is an integral
kernel satisfying the assumptions given in Section \ref{sss.Int_ker}, the above equation
defines the joint intensities. Then, there exists a unique pp $\Phi^{det}$ on $\mR^d$,
such that for each compact $D$, the restriction of $\Phi^{det}$ to $D$
is a determinantal  pp with kernel $K_D$ if and only if the
eigenvalues of $\mathcal{K}_{D}$ are in $[0,1]$. 
\remove{
This latter condition
is equivalent to $\lambda_j^D\in[0,1]$ for all compact $D$;
cf.~\cite[Theorem~4.5.5]{Ben09}. We will call this pp {\em
 determinantal pp with the trace-class integral kernel $K(x,y)$}.
}

\paragraph{Permanental pp}
Similar to the determinantal pp, one  says that a simple pp is a {\em permanental pp} with a kernel $K(x,y)$ with respect to a Radon measure $\mu$ on $\mR^d$ if the  joint intensities of the pp with respect to $\mu^{\otimes k}$
satisfy $\rho^{(k)}(x_1,\ldots,x_k) = \text{per}\big( K(x_i,x_j) \big)_{1
\leq i,j \leq k}$ for all $k$, where $\text{per}\big( \cdot \big)$
stands for the permanent of a matrix. Note that the mean measure of the permanental pp is also  equal to $\alpha(\cdot)=\int_{.}K(x,x)\,\mu(dx)$.
Again, will assume that $K(x,y)$ is an integral
kernel. 
Then, there exists a unique pp $\Phi^{perm}$ on $\mR^d$, such that for each compact $D$, the restriction of $\Phi^{perm}$ to $D$ is a permanental pp with kernel $K_D$; cf.~\cite[Corollary~4.9.9]{Ben09}.  We will call this pp {\em
 permanental pp with the trace-class integral kernel $K(x,y)$}. From~\cite[Proposition~35 and Remark~36]{Ben06}, we also know that
$\Phi^{perm}$ is a Cox pp.\remove{ with intensity field $|F|^2,$ where $F$ is some complex Gaussian process on $\mR^d$.}

\subsection{Comparison results}
The following properties hold true for determinantal and permanental
pp with a trace-class integral kernel $K(x,y)$.

\begin{prop}
\label{p:alpha-weak}
$\Phi^{det}$ is  $\alpha$-weakly sub-Poisson, while
$\Phi^{perm}$ is $\alpha$-weakly super-Poisson;
both comparable with respect to the Poisson pp with mean measure
$\alpha(\cdot)$ given by $\alpha(D)=\int_D
K_D(x,x)\, \mu(\md x)=\sum_{j}\lambda_j^D$, where
the summation is taken over all the eigenvalues $\lambda_j^D$ of
$\mathcal{K}_D$.
\end{prop}

\begin{proof}
\begin{sloppypar}
Since  $K_D(x,y)$ is Hermitian and  positive semi-definite, by Hadamard's
inequality, $\det \big( K_D(x_i,x_j) \big)_{1 \leq i,j \leq k}
\leq \prod_{i=1}^k K_D(x_i,x_i)$ which
implies~(\ref{e.alpha-weakly-sub-Poisson}). For $\Phi^{perm}$, the
proof follows from the permanent analogue of the Hadamard's
inequality (see~\cite{Marvin64}). \qed
\end{sloppypar}
\end{proof}

\begin{prop}
\label{p:nu-weak}
$\Phi^{det}$  is $\nu$-weakly sub-Poisson, while $\Phi^{perm}$
is $\nu$-weakly super-Poisson.
\end{prop}
\begin{proof}
It is known that for each compact $D$, $\Phi^{det}(D) \stackrel{d}{=} \sum_j Bin(1,\lambda_j^D)$ and 
$\Phi^{perm}(D) \stackrel{d}{=} \sum_j Geo(1/(1+\lambda_j^D))$ where the summation is 
taken over all eigenvalues $\lambda_j^D$ of $\mathcal{K}_D$ and 
$Bin(1,\lambda_j^D)$ are independent Bernoulli random variables while
$Geo(1/(1+\lambda_j^D))$ are independent geometric random variables;
cf.~\cite[Theorems~4.5.3 and 4.9.4]{Ben09}. Consequently
$$ \Phi^{det}(D)\le_{cx}Poi(\sum_j\lambda_j^D)\le_{cx}\Phi^{perm}(D)\,, $$
with the left inequality holding provided $\Phi^{det}$ exists (i.e.;
$\lambda_j^D\in[0,1]$ for all compact $D$).
Noting that convex order of integer-valued random variables 
implies ordering of probabilities of taking value~0 concludes the proof (see Proof of \ref{Prop:voids_pp}). \qed
\end{proof}

Alternatively, one can prove the above result via Proposition
\ref{p.association} as the inequality  (\ref{eq:2_void-inequality}) has
been proved for determinantal pp in \cite[Cor. 3.3.]{GeorgiiYoo}.  

\begin{cor} Combining results of Proposition~\ref{p:alpha-weak}
  and~\ref{p:nu-weak} we conclude that
 $\Phi^{det}$ {\em is weakly sub-Poisson}, while $\Phi^{perm}$
{\em is weakly super-Poisson}.
\end{cor}

In the next result, we will
strengthen the above corollary, proving $dcx$ ordering of
finite-dimensional distributions of $\Phi^{det}$ and $\Phi^{perm}$ on mutually disjoint {\em simultaneously observable} sets $D_1,\ldots,D_k$. Simultaneous observability means that the eigenfunctions of $\mathcal{K}_{\bigcup D_i}$, restricted to $D_i$ are also eigenfunctions of $\mathcal{K}_{D_i}$ for every $ i=1,\ldots,k$.

\begin{prop}\label{prop:det-perm-dcx}
Let $\Phi^{det}$ and $\Phi^{per}$ be, respectively,  the
determinantal and permanental pp  with a trace-class integral kernel $K$ and
with $\Phi^{det}$ being defined only if the
spectrum of $\mathcal{K}_{\mR^d}$ is in $[0,1]$.
Denote by $\Phi^{Poi}$ the Poisson pp of mean measure
$\alpha(\cdot)$ given by $\alpha(D)=\sum_{j}\lambda_j^D$ for all
compact $D$, where
the summation is taken over all eigenvalues $\lambda_j^D$ of
$\mathcal{K}_D$. Let $D_1,\ldots,D_k$ be mutually disjoint,
simultaneously observable (with respect to the kernel $K$)
compact subsets of $\mR^d$  and $D = \bigcup D_i$.
Then
\begin{eqnarray*}
\Bigl(\Phi^{det}(D_1),\ldots,\Phi^{det}(D_k)\Bigr)
&\leq_{dcx}& \Bigl(\Phi^{Poi}(D_1),\ldots,\Phi^{Poi} (D_k)\Bigr)\\
&\le_{dcx}&
\Bigl(\Phi^{per}(D_1),\ldots,\Phi^{per}(D_k)\Bigr)\,.
\end{eqnarray*}
\end{prop}

\begin{proof} Let $\{\lambda_{j,i}\}_{j=1,\ldots,J}$  denote the eigenvalues of $\mathcal{K}_{D_i}$ and $\lambda^D_j=\sum_{i=1}^k\lambda^D_{j,i}$, $j=1,\ldots, J$ are the eigenvalues of $\mathcal{K}_D$ with $J$ denoting the number of eigenvalues of $\mathcal{K}_{D}$ ($J=\infty$ and $0$ allowed and in the latter case the sum is understood as $0$).
From~\cite[Prop.~4.5.9]{Ben09}, we know that
$ \Bigl(\Phi^{det}(D_1),\ldots,\Phi^{det}(D_k)\Bigr) \stackrel{d}{=} \sum_{j=1}^J \sum_{i=1}^{N_j} \boldsymbol{\xi}'_{i,j}$,
where $N_j \sim Bin(1,\lambda^D_j)$ and given $N_j$'s, $\boldsymbol{\xi}'_{i,j}$, $i \geq 1$ are independent {\em
 multinomial} vectors $Mul(1,\lambda^D_{j,1}/\lambda^D_j,\ldots,\lambda^D_{j,k}/\lambda^D_j)$~\footnote{
$Mul(k,p_1,\ldots,p_k)$ with $0\le p_i\le 1$, $\sum_{i=1}^kp_i=1$,
has  probability mass  function
$p_{Mul(n,p_1,\ldots,p_k)}(n_1,\ldots,n_k)=\frac{n!}{n_1!\dots n_k!}
p_1^{n_1}\dots p_k^{n_k}$ for $n_1+\ldots+n_k=n$ and 0 otherwise.};  It is easy to see that
$ \Bigl(\Phi^{Poi}(D_1),\ldots,\Phi^{Poi}(D_k)\Bigr) \stackrel{d}{=} \sum_{j=1}^J \sum_{i=1}^{M_j} \boldsymbol{\xi}'_{i,j}$,
where  $M_j \sim Poi(\lambda^D_j)$ with $\boldsymbol{\xi}'_{i,j}$'s and $\lambda^D_{i,j}$'s as defined above.
Due to the independence of $\boldsymbol{\xi_{i,j}}'$'s and the assumption
$\sum_{j=1}^J\lambda_j<\infty$ (local trace-class property of $\mathcal{K}_D$),
it is enough to prove for each $j$, that
$ \sum_{i=1}^{N_j} \boldsymbol{\xi}'_{i,j} \le_{dcx}  \sum_{i=1}^{M_j} \boldsymbol{\xi}'_{i,j}$. Define $g(n) := \EXP{f(\text{sgn}(n) \sum_{i=1}^{|n|} \boldsymbol{\xi}'_{i,j})}$ for $n \in \mZ$ and $f$ $dcx$ function. From Lemmas \ref{lem:MeesterSh-like} and \ref{lem:discrete-convex}, we know that $g(.)$ can be extended to a convex function on $\mR.$ Since we know from (\ref{eqn:sub_Poisson_rvs}) that $Bin(1,\lambda^D_j) \leq_{cx} Poi(\lambda^D_j)$ and hence it follows that $\EXP{g(N_j)} \leq \EXP{g(M_j)}$ as required. This completes the proof of the inequality for the determinantal pp.

Regarding the permanental pp, 
$ \Bigl(\Phi^{per}(D_1),\ldots,\Phi^{per}(D_k)\Bigr) \stackrel{d}{=}
\sum_{j=1}^J  \sum_{i=1}^{K_j} \boldsymbol{\xi}'_{i,j}$ ,
where $K_j \sim Geo(1/(1+\lambda^D_j)$ and $\boldsymbol{\xi}'_{i,j}$'s are as defined 
above; see~\cite[Theorem~4.9.7]{Ben09}.
Similar to the above proof, the required inequality
follows from the ordering, $\sum_{i=1}^{M_j} \boldsymbol{\xi}'_{i,j} \le_{dcx}  \sum_{i=1}^{K_j} \boldsymbol{\xi}'_{i,j}$ for
all $j$, which follows from the fact that
$Poi(\lambda)\le_{cx}Geo(1/(1+\lambda)$ (see
(\ref{eqn:sup_Poisson_rvs})) 
and Lemmas~\ref{lem:MeesterSh-like},~\ref{lem:discrete-convex}.
This completes the proof. \qed 
\end{proof}

\begin{rem}
The key observation used in the above proof
was that the number of points in disjoint,
simultaneously observable sets can be
represented as a sum of independent vectors, which themselves are
binomial (for determinantal) or Poisson (for Poisson) or geometric (for permanental) sums of some further independent vectors. This is exactly the same
representation as for the perturbed pp of Section~\ref{ssec:pl}
(available for any disjoint sets); cf the  proof of
Proposition~\ref{p.pert-lattice}.  In both cases,
this representation and Lemmas~\ref{lem:MeesterSh-like}, \ref{lem:discrete-convex} allow us to conclude $dcx$ ordering of the corresponding vectors.
\end{rem}

\subsubsection*{Example of the  Ginibre process.}
Let $\Phi^{G}$ be the determinantal pp on $\mR^2$
with kernel
$K((x_1,x_2),(y_1,y_2))=\exp[(x_1y_1+x_2y_2)+i(x_2y_1-x_1y_2)]$,
$x_j,y_j\in\mR$,  $j=1,2$,
with respect to the measure
$\mu(d(x_1,x_2))=\pi^{-1}\exp[-x_1^2-x_2^2]\,dx_1dx_2$.
This process is known as the {\em infinite Ginibre} pp.
It is an important example of the determinantal pp
recently studied on the theoretical ground (cf e.g.~\cite{Goldman2010})
and considered in modeling applications (cf.~\cite{Miyoshi2012}). 
Denote by $\Psi^G=\{|X_i|^2: X_i\in\Phi^G\},$ the pp on $\mR^+$
of the squared radii of the points of $\Phi^G$.
This process has an interesting representation  in terms of exponential
random variables, similar to but different from
this of a homogeneous one-dimensional Poisson pp $\Phi^1$;
see~\cite[Theorem~8 (Kostlan)]{Goldman2010}. An interesting questions
posed in~\cite{snorder}
is whether these two processes are $dcx$ ordered. A partial result
given in the cited paper is that 
$\Psi^G([0,r])\le_{cx}\Phi^1([0,r])$ for all $r\ge 0$, was
proved in~\cite{snorder}. 
Full $dcx$ ordering of these two point processes is possible 
studying the simultaneously observable sets for the
Ginibre process.

\begin{cor}
The process of the squared radii of the Ginibre process is
sub-Poisson; i.e.,  $\Psi^G\le_{dcx}\Phi^1$.
\end{cor}
\begin{proof}
We know that an arbitrary finite collection of the annuli centered
at the origin $D_i=\{(x_1,x_2): r_i\le x_1^2+x_2^2\le R_i\}$
is simultaneously observable for this pp;
cf.~\cite[Example~4.5.8]{Ben09}. Using this observation,
Proposition~\ref{prop:det-perm-dcx} and the fact that $dcx$ order
of pp on $\mR$ is generated by the semi-ring of intervals,
we conclude that $\Psi^G$ is $dcx$ smaller than the Poisson pp $\Phi^1$
of unit intensity on $\mR^+$. 
\end{proof} 

\section{Applications and further research}
\label{s.Applications}
In what follows will give some motivating results and preview further ones 
motivating the ideas presented in this paper.

\subsection{Continuum percolation}

\begin{figure}[!t]
\begin{center}
\begin{minipage}[b]{0.30\linewidth}
\caption{\label{f.TwoComponents.subsuperP}
Mean fractions of nodes in the two largest components
of the Boolean models   
generated by perturbed-lattice pp in $d = 2$ with  $Bin(n,1/n)$ 
and $N\!Bin(n,1/(1+n))$ as 
replication kernels, having fixed spherical grains of radius $r$.
The replication kernels converge in $n$, from below and from above in
$dcx$,  respectively,  
to Poisson pp whose  
critical radius is depicted by the dashed line.}
\vspace{8ex}
\end{minipage}
\begin{minipage}[t]{0.6\linewidth}
\includegraphics[width=1\linewidth]{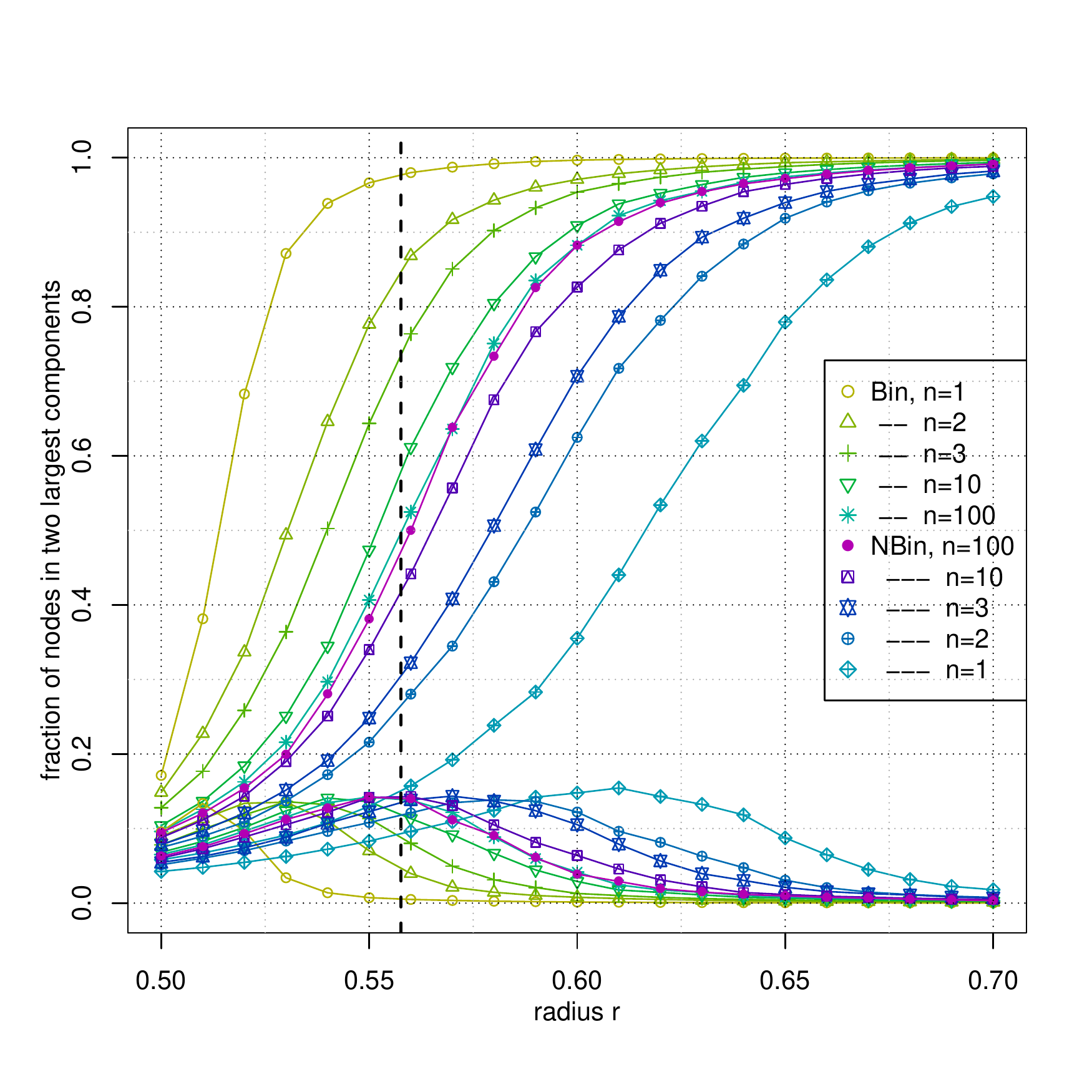}
\end{minipage}
\end{center}
\vspace{-5ex}
  \end{figure}

The Boolean model on a pp $\Phi$ with radius $r$ is defined as $C(\Phi,r) := \bigcup_{X \in \Phi} B_X(r)$, where $B_X(r)$ denotes the ball of radius $r$
centred at $X$. By percolation, we mean the existence of an unbounded
connected subset of the Boolean model. The critical radius for
percolation is defined as $r_c(\Phi) := \inf \{r :
\pr{\mbox{$C(\Phi,r)$ percolates}} > 0 \}$.
We mentioned in the Introduction a heuristic saying that 
clustering worsens percolation. Now, we can use some 
family of perturbed-lattice pp (cf. Section~\ref{ss.pert-Examples}),
 monotone in $dcx$ order, to illustrate this heuristic.
Indeed, Figure~\ref{f.TwoComponents.subsuperP} hints
at ordering of the critical radii of $dcx$ ordered pp in $d = 2$.
However, as shown in~\cite{dcx-perc}, this conjecture is not true in
general: there exists a super-Poisson pp with the critical radius
equal to~0.  What was also shown there, is  that weakly  sub-Poisson pp 
exhibit a (uniformly) non-trivial phase transition in their continuum
percolation model (i.e., admit uniformly non-degenerate lower and
upper bounds for the critical radius). Similar results regarding
$k$-percolation and  SINR-percolation models (arising in 
modeling of connectivity of wireless networks)   
hold for $dcx$-sub Poisson pp.

In what follows, we will present some intuitions leading to the
above results and motivating our special focus on moment
measures and void probabilities in the previous sections.
Specifically, we will introduce 
two newer critical radii $\underline r_c, \overline r_c$,
which act as lower and upper bounds for the usual
critical radius: $\underline r_c\le r_c\le \overline r_c$.
We will show that clustering acts differently on these new radii:
$$\underline r_c(\Phi_2)\le \underline r_c(\Phi_1)\le r_c(\Phi_1)
\le \overline r_c(\Phi_1)\le \overline r_c(\Phi_2)$$
for $\Phi_1$ having smaller voids and moment measures than $\Phi_2$.
This sandwich inequality tels us that $\Phi_1$ exhibits the usual 
phase transition $0<r_c(\Phi_1)<\infty$, provided $\Phi_2$
satisfies a stronger  condition $0<\underline r_c(\Phi_2)$ and
$\overline r_c(\Phi_2)<\infty$. Conjecturing that it holds for Poisson pp
$\Phi_2$, one obtains the  result  on
(uniformly) non-trivial phase transition  for all  weakly sub-Poisson~$\Phi_1$
--- the one proved in~\cite{dcx-perc} in a slightly different way.

\subsubsection{Moment measures and percolation} 
 Let $W_m = [-m,m]^d$ and define $h_{m,k} : (\mR^d)^k \to \{0,1\}$ 
to be the indicator of the event that ${x_1,\ldots,x_k} \in (\Phi \cap
W_m)^k, |x_1| \leq r, \inf_{x \in \partial W_m}|x-x_k| \leq r,
|x_{i+1} - x_i| \leq r \, \, \forall 1 \leq i \leq (k-1)$, where
$\partial W_n$ denotes the boundary of set $W_n$.  
Let $N_{m,k}(\Phi,r) = \sum^{\neq}_{X_1,\ldots,X_k \in \mR^d}h_{m,k}(X_1,\ldots,X_k)$ denote the number of distinct self-avoiding paths of length $k$ from 
the origin $O\in\mR^2$ to the boundary of the box $W_m$ in the Boolean model and $N_m(\Phi,r) = \sum_{k \geq 1} N_{m,k}(\Phi,r)$ to be the total number of distinct self-avoiding paths 
to the boundary of the box. We define the following ``lower'' critical radius:
$$\ur_c(\Phi)  :=  \inf \{r : \liminf_m \EXP{N_m(\Phi,r)} > 0 \}\,.$$
Note that $r_c(\Phi)=\inf \{r : \lim_m \pr{N_m(\Phi,r) \geq 1}
> 0 \}$, with the limit existing because the events $\{N_m(\Phi,r) \geq 1 \}$
form a decreasing sequence in $m$, and by Markov's inequality, we have that 
indeed $\ur_c(\Phi) \leq r_c(\Phi)$ for a stationary pp $\Phi$.

\begin{prop}
\label{prop:lower_crit_rad}
Let  $C_j=C(\Phi_j,r)$, $j=1,2$ be two Boolean models with simple pp of germs
$\Phi_j$, $j=1,2,$ and $\sigma$-finite
$k\,$th moment measures $\alpha_j^{k}$ for all $k\ge1$ respectively.
If $\alpha_1^{(k)}(\cdot) \leq \alpha_2^{(k)}(\cdot)$ 
for all $k\ge 1$, then $\ur_c(\Phi_1)\ge \ur_c(\Phi_2).$
In particular, for a stationary, $\alpha$-weakly sub-Poisson pp $\Phi_1$ of unit
intensity we have that
$\theta_d \ur_c(\Phi_1)^d \geq 1$ where $\theta_d$ is the volume of the unit ball.~\footnote{Similar
to open paths from the origin   to $\partial Q_m$, one can define an {\em open path on the germs of $\Phi$
crossing the rectangle  $ [0,m] \times [0,3m] \times \ldots \times [0,3m]$ across the shortest side} and
define yet another critical radius $r_s(\Phi)$ as the smallest $r$ for which such a path
exists with positive probability for an arbitrarily large $m$ (\cite[(3.20)]{MeeRoy96}). An analogous inequality holds true for this critical radius too.} 
\end{prop}
\begin{proof}
The proof relies on the following easy derivation of the closed form
expressions for $\EXP{N_m(\Phi_j,r)}, j=1,2$: $\EXP{N_m(\Phi,r)} =
\sum_{k \geq 1} \EXP{N_{m,k}(\Phi,r)}$ and 
$$\EXP{N_{m,k}(\Phi,r)}=\int_{(\mR^d)^k}h_{m,k}(x_1,\ldots,x_k) \alpha_j^{(k)}(dx_1,\ldots,dx_k)\,.$$ 
For the second part of the proof note that the above summation over $k$ can be
taken over $k \geq m_r:=\lfloor m/r \rfloor - 1$, where $\lfloor
a\rfloor$  denotes the larges integer not larger than $a$. 
Indeed, the maximal distance that can be reached by a path of length
$k$ in $C(\Phi,r)$ is $(k+1)r$ and hence $h_{m,k} \ge 1$ implies that $k
\geq m_r$.  Consequently, for $\alpha$-weakly sub-Poisson pp $\Phi$
\begin{eqnarray*}
\EXP{N_m(\Phi,r)} & \leq & \sum_{k \geq m_r}
\int_{(\mR^d)^k}h_{m,k}(x_1,\ldots,x_k) dx_1\ldots dx_k \\
&\leq&  \sum_{k \geq m_r} (\theta_d r^d)^k= \frac{(\theta_d
  r^d)^{m_r}}{1 - \theta_dr^d}\,,
\end{eqnarray*}
where the second inequality  follows by releasing the condition that
$x_k$ is close to $\partial W_m$.
Thus, $\EXP{N_m(\Phi,r)}<\infty$ for $\theta_d r^d < 1$ and hence 
the result  $\theta_d \ur_c(\Phi)^d \geq 1$. \qed
\end{proof}

An interesting consequence of the above result is that $r_c(\Phi) \geq \theta_d^{-\frac{1}{d}} \to \infty$ as $d \to \infty$ 
for $\alpha$-weakly sub-Poisson pp whereas $r_c(\mZ^d) = \frac{1}{2}$ for all $d \geq 1$ i.e, achieving percolation on a sub-Poisson pp
is distinctly more difficult than on a regular lattice in higher
dimensions. This was already known for Poisson pp
(see~\cite{Penrose96})
and now it shows
that the $\alpha$-weakly sub-Poissonianity  does not help in (prevents
from!) percolating faster.   

Given a graph, let $c_n(G)$ be the expected number of self-avoiding walks starting from a fixed point in the lattice. Then the expected {\em connective constant} of the graph is $\mu(G) := \lim_n c_n(G)^{\frac{1}{n}}$. From the proof above, one can also infer that $c_n(C(\Phi_1,r)) \leq c_n(C(\Phi_2,r))$ for $\alpha_1^{(n)}(\cdot) \leq \alpha_2^{(n)}(\cdot)$ and $\mu(C(\Phi,r)) \leq \theta_dr^d$ for a $\alpha$-weakly sub-Poisson pp $\Phi$. 

\subsubsection{Void probabilities and percolation}
Though we are interested in the percolation of Boolean models (continuum percolation models),
but as is the wont in the subject we shall use discrete
percolation models as approximations. For $r > 0, x \in \mR^d$, define the following subsets of $\mR^d$ : 
$Q^r := (-\frac{1}{2r},\frac{1}{2r}]^d$ and  $Q^r(x) := x + Q^r.$
We will consider the following  discrete graph parametrized by  $n \in \mN$ :
$\mL^{*d}_n=(\mZ^d_n, \mE^{*d}_n)$ is the usual  {\em close-packed lattice graph} scaled down by
the factor $1/n$. It has $\mZ^d_n  =\frac{1}{n}\mZ^d$, where $\mZ$ is
the set of integers, as the set of vertices  and the set of edges $\mE^{*d}_n := \{ \langle z_i,z_j \rangle \in
(\mZ^d_n)^2 :  Q^{\frac{n}{2}}(z_i) \cap Q^{\frac{n}{2}}(z_j) \neq \emptyset
\}$.

A {\em contour} in $\mL^{*d}_n$ is a minimal collection of vertices such
that any infinite path in $\mL^{*d}_n$
from the origin has to contain one of these
vertices (the minimality condition implies that the removal of any vertex
from the collection will lead to existence of an infinite path from
the origin without any intersection with the remaining vertices in the
collection). Let $\Gamma_n$ be the set of all contours around the origin in
$\mL^{*d}_n$. For any subset of points $\gamma\subset\mR^d$, in particular for paths
$\gamma \in \Upsilon^n_m(K),\, \Gamma_n$, we define  $Q_{\gamma} = \bigcup_{z \in \gamma} Q^n(z)$.

With these notations, we can define the ``upper'' critical radius $\ovr_c(\Phi)$.
\begin{equation}
\label{e:urc}
\ovr_c=\ovr_c(\Phi):= \inf \Bigl\{r>0 : \text{for all}\;  n\ge1, \sum_{\gamma \in
\Gamma_n}\pr{C(\Phi,r) \cap Q_{\gamma} = \emptyset} <  \infty\Bigr\}\,.
\end{equation}
It might be seen as the critical radius corresponding to the phase transition when the discrete model $\mL^{*d}_n=(\mZ^d_n, \mE^{*d}_n)$,  approximating
$C(\Phi,r)$ with an arbitrary precision, starts percolating through the Peierls argument. As a consequence, $\ovr_c(\Phi) \geq r_c(\Phi)$ (see \cite[Lemma 4.1]{perc-dcx}).
The following ordering result follows immediately from the definition.
\begin{cor}
\label{cor:ord_crit_rad_racs}
Let  $C_j=C(\Phi_j,r)$, $j=1,2$ be two Boolean models with simple pp of germs
$\Phi_j$, $j=1,2$. 
If $\Phi_1$ has smaller voids probabilities than $\Phi_2$
then $\ovr_c(\Phi_1)\le\ovr_c(\Phi_2)$.
\end{cor}

\begin{rem}
Even if the finiteness of 
$\ovr_c$ 
is not clear even for Poisson pp  and hence
Corollary~\ref{cor:ord_crit_rad_racs}
cannot be directly used to prove the
finiteness of the critical radii of $\nu$-weakly sub-Poisson pp,
the approach based on void probabilities
can be refined, as shown in~\cite{dcx-perc}, to conclude the
aforementioned  property. 
\end{rem}

\subsection{Multiple  coverage}
For a point process $\Phi$, define the $k$-covered set $C_k(\Phi,r) :=
\{ x \in \mR^d : \exists X_1,\ldots,X_m \in \Phi, m \geq k \ni x \in 
\bigcap_{i=1}^m B_{X_i}(r) \}$. Heuristically, clustering should reduce the $1$-covered region but 
increase the $k$-covered region for large $k$ and we present a more formal statement of the same. Expected volume of the 
$k$-covered region is one of the important quantities of interest in sensor networks and our result has obvious implications
regarding the choice of $k$ or the point process in the context of sensor networks. We introduce another stochastic order to state the result. We say that two random variables
$X,Y$ are ordered in {\em uniformly convex variable order} (UCVO)($X \leq_{uv} Y$) if their respect density 
funtions $f,g$ satisfy the following conditions : $supp(f) \subset supp(g)$, $f(\cdot)/g(\cdot)$ is an unimodal function but their respective distribution functions are not ordered i.e,
$F(\cdot) \nleqslant G(\cdot)$ or vice-versa (see \cite{Whitt1985}) and where $supp(.)$ denotes the support of a function.
Denote by $\|A\|$ the Lebesgue's measure of bBs $A\subset\mR^d$.
\begin{prop}
\label{prop:k_coverage}
Let $\Phi_1$ and $\Phi_2$ be two simple, stationary pp such that
$\Phi_1(B_O(r)) \leq_{uv} \Phi_2(B_O(r)$ for $r \geq 0$.
Then there exists $k_0 \geq 1$ such that for any 
bBs $W \subset \mR^d$
\begin{eqnarray*}
\forall_{k: 1 \leq k \leq k_0}, & & \EXP{\|C_k(\Phi_1,r) \cap W \|} \geq \EXP{\|C_k(\Phi_2,r) \cap W \|} \, \, \mbox{and} \\
\forall_{k > k_0}, & & \EXP{\|C_k(\Phi_1,r) \cap W \|} \leq \EXP{\|C_k(\Phi_2,r) \cap W \|}.
\end{eqnarray*}
\end{prop}
\begin{proof}
\begin{sloppypar}
Firstly note that
$\EXP{\|C_k(\Phi_i,r) \cap W \|} = \int_{W} \pr{\Phi_i(B_x(r)) \geq k}
dx = \|W\| \pr{\Phi_i(B_O(r)) \geq k}$  for $i=1,2$.
Now it suffices to show that $\pr{\Phi_1(B_O(r)) \geq k} - \pr{\Phi_2(B_O(r)) \geq k}$ changes sign exactly once in $k$ for $k \geq 1$. This is implied by the UCVO order 
(see \cite[Section 2 and Theorem 1]{Whitt1985}).
\end{sloppypar}
\end{proof}
It is known that log-concavity of $f/g$ implies UCVO order as well as convex ordering. We have used the latter implication in our examples for sub-Poisson (see Section \ref{sec:sub_poisson_lattice}) and 
super-Poisson perturbed lattices (see Section \ref{sec:sup_poisson_lattice}). We can take for $\Phi_1$ any of the sub-Poisson perturbed lattices presented in this article or 
determinantal pp and $\Phi_2$ as a Poisson pp. We can also take $\Phi_1$ to be a Poisson pp and $\Phi_2$ to be any of the super-Poisson perturbed lattices presented in this article or 
permenantal pp. 

\subsection{Further applications}
\subsubsection{Minimal spanning forest (MSF)}
In a recent work~\cite{Hirsch2011}, the authors show 
the connectivity of some approximations of the MSF
for the weakly sub-Poisson point processes (cf. 
the  conjecture by Aldous and Steele (\cite{Aldous92}) that the MSF of Poisson pp
is almost surely connected, proved  by Alexander (\cite{Alexander95}) for dimension $d=2$
in 1995).

\subsubsection{First passage percolation}
Existence of arbitrarily large voids in Poisson pp
was shown in~\cite{BBMfppSINR} to be a reason 
of infinite end-to-end packet-delivery delays in a 
time-space  SINR model, studied in the framework of first passage
percolation problem. Superposing the Poisson pp with 
an independent lattice of arbitrarily small intensity 
makes the delays finite.  The latter result remains true 
when using  a simple perturbed lattice, in which case the
superposition is an example of a ($dcx$) sub-Poisson pp. Generalization to 
an arbitrary sub-Poisson pp is an open question. An interesting connection exists to the work of \cite{Vandenberg93,Marchand02} on inequalities for time constants in first passage percolation on $\mZ^d$ with differing 
distributions for edge-passage times. More precisely, it was shown that more variable (in the sense of convex order) edge-passage times lead to faster transmission i.e, smaller time constant. An analogous result for time constants in the continuum case would be a welcome addition to the subject.

\subsubsection{Lilypond growth model}
In \cite[Section 4.2]{Daley05}, it is shown that the Lilypond growth model exists for sub-Poisson pp, though not using the same terminology. Our examples of sub-Poisson pp adds to the list
of examples given in \cite{Daley05} for which Lilypond growth model exists. Further, it was shown that sub-Poisson pp with absolutely continuous (w.r.t. Lebesgue measure) $\alpha^{(k)}(.)$'s do not percolate.
The absolute continuity condition also holds true for our examples.   

\subsubsection{Random geometric complexes}
This topological extension of random geometric graphs (\cite{Penrose03}) was introduced and studied in \cite{Kahle11} on Poisson pp exploiting the connection between the Betti numbers of a random geometric complex and component counts of the corresponding random geometric graph (\cite[Chapter 3]{Penrose03}). The motivation lies in the recent subject of topological data analysis. In an upcoming work (\cite{Adler12}), we study these models on more general stationary point processes using tools of stochastic ordering as well as asymptotic analysis of joint intensities and void probabilities. In particular, if we denote $r^{con}_n(\Phi)$ as the critical contractibility radius for the \u{C}ech complex on $\Phi \cap [-\frac{n^{\frac{1}{d}}}{2},\frac{n^{\frac{1}{d}}}{2}]$ (i.e, the least radius $r$ above which the Boolean model $C(\Phi \cap [-\frac{n^{\frac{1}{d}}}{2},\frac{n^{\frac{1}{d}}}{2}],r)$ becomes homotopic to a single point), then $r^{con}_n(\Phi) =O((\log n)^{\frac{1}{d}})$ for a $\nu$-weakly sub-Poisson 
pp whereas the critical contractibility radius of a Poisson pp is $\Theta((\log n)^{\frac{1}{d}})$. For \u{C}ech and Vietoris-Rips complexes on $\alpha$-weakly sub-Poisson pp, it is shown that order of the radii for existence of non-zero $k$th Betti numbers ($k \geq 1$) are $\Omega(n^{-\frac{1}{d(k+1)}})$ and $\Omega(n^{-\frac{1}{d(2k+1)}})$ respectively i.e, at least that of the Poisson pp. For specific weak sub-Poisson pp such as the Ginibre determinantal pp for which one has more accurate information about its joint intensities and void probabilities, it is shown that the correct orders differ significantly from that of the Poisson pp. The stronger assumption of negative association allows one to obtain variance bounds and hence derive asymptotics for existence of Betti numbers with high probability in the intermediate regime.  

Subgraph counts and connective constants of random geometric graphs constitute specific instances of order statistics of pps. Scaling limits of order statistics of Poisson pp has garnered some interest in recent times ; see \cite{Schulte12}. Our techniques can easily yield that first moments of the order statistics are ordered for $\alpha$-weakly ordered pp but the question of further asymptotics remains open.

\subsubsection{Applications in modeling}
In the context of wireless networks,
pp are used to model locations of emitters/receivers.
An ubiquitous assumption when modeling 
base stations in cellular networks
is to consider deterministic lattices (usually hexagonal). 
On the other hand, mobile users 
are usually modeled by a Poisson pp. 
Both the assumptions are too simplistic.
In reality, patterns of base stations are neither perfectly periodic,
due to various locational constants nor completely
independent because of various interactions:
social, human interactions typically introduce more clustering, 
while the medium access protocols implemented in mobile wireless 
devices (as e.g. CSMA used in the popular WiFi technology)
tend to separate active users. 
One clearly sees the interest in  perturbed-lattice models in
this context. We believe also that  our work may lay the groundwork 
in other domains, e.g. in social and economic sciences,
where one studies the impact of 
clustering on the macroscopic properties of
models (cf. e.g.~\cite{CL11}).

\subsubsection{Further research}
Another motivation to study sub-Poisson perturbed lattices comes
from their relations to zeros of {\em Gaussian analytic functions} (GAF), 
cf~\cite{Sodin06}, 
whose points exhibit repulsion at smaller distances and independence over large distances.
However, the points seem more regularly distributed than in 
Poisson pp (\cite{Peres05}). This asks the question whether zeros of GAF 
are comparable in some sense to  Poisson pp. 
{\em Gibbsian pp} is another well-known class of point processes,
which depending on the nature of the potential would be more or
less clustering. Super and sub-poissonianity (even in the weak sense) have
not been studied yet for Gibbsian pp. Devising {\em statistical tests} for sub-Poissonianity would be desirable.

\appendix
\section*{Appendix}
\setcounter{section}{1}
\pdfbookmark[1]{Appendix}{appendix} 
\vspace{-0.1in}
The following result, similar to~\cite[Lem\-ma~2.17]{MeesterSh93}
is used in the proof of Propositions~\ref{p.pert-lattice}
and~\ref{prop:det-perm-dcx}.
\begin{lem}
\label{lem:MeesterSh-like}
Let ${\boldsymbol\xi}_i=(\xi_i^1,\ldots,\xi_i^k)\in\mR^k$, ($i \in \mZ$)
be independent, identically distributed  vectors of (possibly
dependent) non-negative random variables.
Suppose $f$ is a $dcx$ function on $\mR^k$. Then, the function $g$ defined
on $\mZ$ by
 $g(n) =\EXP{f(\text{sgn}(n)\sum_{i=1}^{|n|}{\boldsymbol\xi}_i)}$ for $n\not=0$ and  $g(0)=0$ is convex on~$\mZ$.
\end{lem}
\begin{proof}
We will prove that $g(n)$ has non-negative
second differences
\begin{equation}\label{eq:second-differences}
g(n-1)+g(n+1)-2g(n)\ge 0\qquad \text{for all $n\in\mZ$}
\end{equation}
and use the first part of Lemma~\ref{lem:discrete-convex}.
To prove~(\ref{eq:second-differences}),
define $G(n,m) := \sum_{i=n+1}^m{\boldsymbol\xi}_i$
for $0\le n < m$ and $G(n,n):=(0,0,\ldots,0)\in\mR^k$ for $n\ge0$.
We have for $n \ge 1$,
\begin{eqnarray*}
2g(n)&=&2\EXP{f\Bigl(G(0,n)\Bigr)}\\
&=&\EXP{f\Bigl(G(0,n-1)+G(n-1,n)\Bigr)}+
\EXP{f\Bigl(G(0,n)\Bigr)}\\
&=&\EXP{f\Bigl(G(0,n-1)+G(n,n+1\Bigr)}+
\EXP{f\Bigl(G(0,n)\Bigr)}\\
&=&\EXP{f\Bigl(G(0,n-1)+G(n,n+1\Bigr)+
f\Bigl(G(0,n)\Bigr)}\\
&\le &\EXP{f\Bigl(G(0,n-1)\Bigr)+
f\Bigl(G(0,n)+G(n,n+1)\Bigr)}\\
&=&g(n-1)+g(n+1)\,,
\end{eqnarray*}
where for the third equality
we have used mutual independence of $G(0,n-1), G(n-1,n),
G(n,n+1)$ and the fact that $G(n-1,n)$ and $G(n,n+1)$ have the same
distribution, while the inequality follows from the
$dcx$ property of $f$ and the assumption
${\boldsymbol\xi}_{i}\ge0$.
This proves~(\ref{eq:second-differences})
for $n\ge1$. Similar reasoning allows to
show~(\ref{eq:second-differences}) for $n\le-1$.
Finally, note that for  $n=0$
\begin{eqnarray*}
2g(0)&=&2f\Bigl((0,\ldots,0)\Bigr)\\
&=&\EXP{f\Bigl(-G(0,1)+G(0,1)\Bigr)+
f\Bigl((0,\ldots,0)\Bigr)}\\
&\le &\EXP{f\Bigl(-G(0,1)\Bigr)+
f\Bigl(G(0,1)\Bigr)}\\
&=&g(-1)+g(1)\,,
\end{eqnarray*}
\end{proof}

We will prove the following two technical results regarding convex
functions. We were not able to find their proofs  in the literature.
\begin{lem}\label{lem:discrete-convex}
Let $g(n)$ be a real valued function defined for all integer
$n\in\mZ$ and satisfying condition~(\ref{eq:second-differences}).
Then 
for all $n\ge 2$
\begin{equation}\label{eq:discrete-convex}
g\Bigl(\sum_{i=1}^n\lambda_i k_i\Bigr)\le
\sum_{i=1}^n\lambda_i g(k_i)
\end{equation}
for all $k_i\in\mZ$ and $0\le \lambda_i\le 1$,
 $\sum_{i=1}^n\lambda_i=1$ such that
 $\sum_{i=1}^n\lambda_i k_i\in\mZ$.
Moreover, function $g(\cdot)$ can be extended
to a real valued convex function defined on real numbers~$\mR$.
 \end{lem}
\begin{proof}
As mentioned in~\cite[Section~V.16.B.10.a]{MOA2009}
it is easy to see that~(\ref{eq:second-differences})
is equivalent to~(\ref{eq:discrete-convex}) with $n=2$.
Assume now that~(\ref{eq:discrete-convex}) holds true for some
$n\ge2$ (and all $0\le\lambda_i\le1$, $k_i\in\mZ$, $i=1,\ldots,n$ satisfying
$\sum_{i=1}^n\lambda_i=1$,  $\sum_{i=1}^n\lambda_i k_i\in\mZ$).
We will prove that it holds true for $n+1$ as well.
In this regard,
define  for a given  $k\in\mZ$ and distinct (otherwise we use directly
the inductive assumption) $k_1,\ldots,k_{n+1}\in\mZ$,
the following functions:
\begin{eqnarray*}
\lambda_{n}=\lambda_{n}(\lambda_1,\ldots,\lambda_{n-1})&:=&
\frac{k-k_{n+1}-\sum_{i=1}^{n-1}\lambda_i(k_i-k_{n+1})}{k_n-k_{n+1}}\\
\lambda_{n+1}=
\lambda_{n+1}(\lambda_1,\ldots,\lambda_{n-1})&:=&1-\sum_{i=1}^{n-1}\lambda_i-
\lambda_{n}(\lambda_1,\ldots,\lambda_{n-1})\\
F(\lambda_1,\ldots,\lambda_{n-1})&:=&\sum_{i=1}^{n-1}\lambda_ig(k_i)
+\lambda_{n}(\lambda_1,\ldots,\lambda_{n-1})g(k_n)\\
&&\hspace{5em}+\lambda_{n+1}(\lambda_1,\ldots,\lambda_{n-1})g(k_{n+1})\,.
\end{eqnarray*}
Note that for any $\lambda_1,\ldots,\lambda_{n-1}$ we have
$\sum_{i=1}^{n+1}\lambda_i=1$ and $\sum_{i=1}^{n+1}\lambda_ik_i=k$.
Consider the following
subset of the $n-1$-dimensional unit cube
$C:=\Bigl\{(\lambda_1,\ldots,\lambda_{n-1})\in[0,1]^{n-1}:
0\le \lambda_{n}\le 1, 0\le \lambda_{n+1}\le 1\Bigr\}$.
The proof of the inductive step will be
completed if we show that
$F(\cdot) \ge g(k)$ on $C$.
In this regard note that $C$ is {\em closed} and {\em  convex}.
Assume moreover that $C$ is not empty; otherwise the
condition~(\ref{eq:discrete-convex})
is trivially satisfied.
Note also that  $F(\cdot)$ is  an affine, real valued function
defined on $\mR^{n-1}$. Hence, by the maximum principle,
the affine (hence convex) function $-F$
attains its maximum relative to $C$ on some point
$(\lambda_1^0,\ldots,\lambda_{n-1}^0)\in\partial C$ of the boundary of
$C$. Consequently, we have
$F(\cdot)\ge F(\lambda_1^0,\ldots,\lambda_{n-1}^0)$ on $C$ and the
 proof of the inductive step will be completed if we show that
$ F(\lambda_1^0,\ldots,\lambda_{n-1}^0)\ge g(k)$.
In this regard, denote
$\lambda_{n}^0=\lambda_n(\lambda_1^0,\ldots,\lambda_{n-1}^0)$ and
$\lambda_{n+1}^0=\lambda_{n+1}(\lambda_1^0,\ldots,\lambda_{n-1}^0)$.
Using the continuity of the functions
$\lambda_{n}(\cdot)$ and $\lambda_{n+1}(\cdot)$
 is not difficult to verify that
 $(\lambda_1^0,\ldots,\lambda_{n-1}^0)\in\partial C$ implies
$\lambda_j^0=0$  for some $j=1,\ldots,n+1$.
Thus, by
our inductive assumption, $F(\lambda_1^0,\ldots,\lambda_{n-1}^0)=
\sum_{i=1, i\not=j}^{n+1}\lambda_i^0g(k_i)\ge g(k)$,
which completes the proof of~(\ref{eq:discrete-convex}) for all $n\ge2$.

For the second statement,
we  recall the arguments used in~\cite{Yan1997}
to show that a function satisfying~(\ref{eq:discrete-convex})
for all $n\ge 2$ (called {\em globally convex function} there)
has a convex extension on $\mR$. In this regard,
consider the epigraph $\text{epi}(g):=\{(k,\mu)\in\mZ\times\mR: \mu\ge
g(k)\}$ of $g$ and its convex envelope $\text{epi}^{co}(g)$.
It is easy to see that
$\text{epi}^{co}(g)=\{(x,\mu)\in\mR^2:
\mu\ge\sum_{i=1}^n\lambda_ig(k_i)$ for some $k_i\in\mZ$, $0\le
\lambda_i\le 1$, $\sum_{i=1}^n\lambda_i=1$ and $\sum_{i=1}^n\lambda_ik_i=x\}$.
Define $\tilde g(x):=\inf\{\mu:(x,\mu)\in\text{epi}^{co}(g)\}$ for all
$x\in\mR$. The convexity of $\text{epi}^{co}(g)$ implies that $\tilde
g$ is convex on $\mR$ and the global
convexity~(\ref{eq:discrete-convex}) of $g$ implies that $\tilde g$ is an
extension of $g$.
This completes the proof.
\end{proof}

\ack
\pdfbookmark[1]{Acknowledgments}{acknowledgments} 

The authors wish to thank Manjunath Krishnapur for introducing them to perturbed lattices and answering various queries on determinantal point processes. Also, the results on negatively associated point processes originated from discussions with Manjunath Krishnapur and Subhrosekhar Ghosh. DY was supported by INRIA  and ENS Paris where most of this work was done. He is also thankful to research grants from EADS(Paris), Israel Science Foundation (No: 853/10) and AFOSR (No: FA8655-11-1-3039). The authors also wish to thank one of the anonymous referees for pointing out various connections of our work to existing literature such as high-dimensional percolation, first-passage percolation and Lilypond growth model.  

\bibliographystyle{apt}
\pdfbookmark[1]{References}{references} 


\end{document}